\documentclass[12pt]{article} 
\usepackage{latexsym} 
\usepackage{amsmath}
\usepackage{mathptmx}
\usepackage{amssymb} 
\usepackage{amssymb} 
\usepackage{amsfonts}
\usepackage{amscd} 
\usepackage[all]{xy}
\usepackage{graphicx}
\usepackage{layout}
\begin{document} 
\newtheorem{Th}{Theorem}[section]
\newtheorem{Cor}{Corollary}[section]
\newtheorem{Prop}{Proposition}[section]
\newtheorem{Lem}{Lemma}[section]
\newtheorem{Def}{Definition}[section]
\newtheorem{Rem}{Remark}[section]
\newtheorem{Ex}{Example}[section]
\newtheorem{stw}{Proposition}[section]


\newcommand{\bet}{\begin{Th}}
\newcommand{\ent}{\stepcounter{Cor}
   \stepcounter{Prop}\stepcounter{Lem}\stepcounter{Def}
   \stepcounter{Rem}\stepcounter{Ex}\end{Th}}


\newcommand{\bec}{\begin{Cor}}
\newcommand{\enc}{\stepcounter{Th}
   \stepcounter{Prop}\stepcounter{Lem}\stepcounter{Def}
   \stepcounter{Rem}\stepcounter{Ex}\end{Cor}}
\newcommand{\bep}{\begin{Prop}}
\newcommand{\enp}{\stepcounter{Th}
   \stepcounter{Cor}\stepcounter{Lem}\stepcounter{Def}
   \stepcounter{Rem}\stepcounter{Ex}\end{Prop}}
\newcommand{\bel}{\begin{Lem}}
\newcommand{\enl}{\stepcounter{Th}
   \stepcounter{Cor}\stepcounter{Prop}\stepcounter{Def}
   \stepcounter{Rem}\stepcounter{Ex}\end{Lem}}
\newcommand{\bef}{\begin{Def}}
\newcommand{\enf}{\stepcounter{Th}
   \stepcounter{Cor}\stepcounter{Prop}\stepcounter{Lem}
   \stepcounter{Rem}\stepcounter{Ex}\end{Def}}
\newcommand{\ber}{\begin{Rem}}
\newcommand{\enr}{
   \stepcounter{Th}\stepcounter{Cor}\stepcounter{Prop}
   \stepcounter{Lem}\stepcounter{Def}\stepcounter{Ex}\end{Rem}}
\newcommand{\bee}{\begin{Ex}}
\newcommand{\ene}{
   \stepcounter{Th}\stepcounter{Cor}\stepcounter{Prop}
   \stepcounter{Lem}\stepcounter{Def}\stepcounter{Rem}\end{Ex}}
\newcommand{\Proof}{\noindent{\it Proof\,}:\ }
\newcommand{\beP}{\Proof}
\newcommand{\enP}{\hfill $\Box$ \par\vspace{5truemm}}

\newcommand{\EE}{\mathbb{E}}
\newcommand{\QQ}{\mathbb{Q}}
\newcommand{\R}{\mathbb{R}}
\newcommand{\C}{\mathbb{C}}
\newcommand{\ZZ}{\mathbb{Z}}
\newcommand{\KK}{\mathbb{K}}
\newcommand{\NN}{\mathbb{N}}
\newcommand{\PP}{\mathbb{P}}
\newcommand{\HH}{\mathbb{H}}
\newcommand{\uuu}{\boldsymbol{u}}
\newcommand{\xxx}{\boldsymbol{x}}
\newcommand{\aaa}{\boldsymbol{a}}
\newcommand{\bbb}{\boldsymbol{b}}
\newcommand{\AAA}{\mathbf{A}}
\newcommand{\BBB}{\mathbf{B}}
\newcommand{\LLL}{\mathbf{L}}
\newcommand{\ccc}{\boldsymbol{c}}
\newcommand{\iii}{\boldsymbol{i}}
\newcommand{\jjj}{\boldsymbol{j}}
\newcommand{\kkk}{\boldsymbol{k}}
\newcommand{\rrr}{\boldsymbol{r}}
\newcommand{\FFF}{\boldsymbol{F}}
\newcommand{\yyy}{\boldsymbol{y}}
\newcommand{\ppp}{\boldsymbol{p}}
\newcommand{\qqq}{\boldsymbol{q}}
\newcommand{\nnn}{\boldsymbol{n}}
\newcommand{\vvv}{\boldsymbol{v}}
\newcommand{\eee}{\boldsymbol{e}}
\newcommand{\fff}{\boldsymbol{f}}
\newcommand{\www}{\boldsymbol{w}}
\newcommand{\0}{\boldsymbol{0}}
\newcommand{\lon}{\longrightarrow}
\newcommand{\ga}{\gamma}
\newcommand{\pa}{\partial}
\newcommand{\QED}{\hfill $\Box$}
\newcommand{\id}{{\mbox {\rm id}}}
\newcommand{\Ker}{{\mbox {\rm Ker}}}
\newcommand{\Image}{{\mbox {\rm Image}}}
\newcommand{\grad}{{\mbox {\rm grad}}}
\newcommand{\ind}{{\mbox {\rm ind}}}
\newcommand{\rot}{{\mbox {\rm rot}}}
\newcommand{\diver}{{\mbox {\rm div}}}
\newcommand{\Gr}{{\mbox {\rm Gr}}}
\newcommand{\GL}{{\mbox {\rm GL}}}
\newcommand{\LG}{{\mbox {\rm LG}}}
\newcommand{\Diff}{{\mbox {\rm Diff}}}
\newcommand{\Symp}{{\mbox {\rm Symp}}}
\newcommand{\Ct}{{\mbox {\rm Ct}}}
\newcommand{\Uns}{{\mbox {\rm Uns}}}
\newcommand{\rank}{{\mbox {\rm rank}}}
\newcommand{\sign}{{\mbox {\rm sign}}}
\newcommand{\Spin}{{\mbox {\rm Spin}}}
\newcommand{\Sp}{{\mbox {\rm Sp}}}
\newcommand{\Int}{{\mbox {\rm Int}}}
\newcommand{\Hom}{{\mbox {\rm Hom}}}
\newcommand{\Tan}{{\mbox {\rm Tan}}}
\newcommand{\codim}{{\mbox {\rm codim}}}
\newcommand{\ord}{{\mbox {\rm ord}}}
\newcommand{\Iso}{{\mbox {\rm Iso}}}
\newcommand{\corank}{{\mbox {\rm corank}}}
\def\mod{{\mbox {\rm mod}}}
\newcommand{\pt}{{\mbox {\rm pt}}}
\newcommand{\qed}{\hfill $\Box$ \par}
\newcommand{\spe}{\vspace{0.4truecm}}
\renewcommand{\0}{\mathbf 0}
\newcommand{\ad}{{\mbox{\rm ad}}}
\newcommand{\xdownarrow}[1]{%
  {\left\downarrow\vbox to #1{}\right.\kern-\nulldelimiterspace}
}

\newcommand{\dint}[2]{{\displaystyle\int}_{{\hspace{-1.9truemm}}{#1}}^{#2}}

\title{Cofrontals
}

\author{Goo \textsc{Ishikawa}\thanks{Faculty of Science, Department of Mathematics, Hokkaido University, Sapporo 060-0810, Japan. e-mail: 
ishikawa@math.sci.hokudai.ac.jp
}}


%

\renewcommand{\thefootnote}{\fnsymbol{footnote}}
\footnotetext{
\\
\qquad
The author was supported by JSPS KAKENHI No.15H03615.}

\date{}

\maketitle

\begin{abstract} 
In this paper we introduce the notion of cofrontal mappings, 
as the dual objects to frontal mappings, 
and study their basic local and global properties.  
Cofrontals are very special mappings 
and far from generic nor stable except for the case of submersions. 
It is observed 
that any smooth mapping can be $C^0$-approximated by a possibly \lq\lq unfair" cofrontal or a frontal. However global \lq\lq fair" cofrontals are very restrictive to exist. 
Then we give a method to construction \lq\lq fair" cofrontals with fiber-dimension one 
and a target-local diffeomorphism classification of such cofrontals, 
under some finiteness condition. 
\end{abstract}

\section{Introduction}

In the previous papers (see \cite{Ishikawa13}\cite{Ishikawa14}\cite{Ishikawa18}\cite{Ishikawa18-2}) 
we have introduced and studied the notion of frontal map-germs. 
A map-germ $f : (N, a) \to (M, b)$ from an $n$-dimensional manifold $N$ to an $m$-dimensional manifold $M$ is called a {\it frontal} if $n \leq m$ and 
there exists a smooth $n$-plane field $\widetilde{f}$ along $f$, 
i.e. which commutes 
$$
\xymatrix{
 & {\mbox{\rm Gr}}(n, TM) \ar[d]^{\pi}
 \\
(N, a) \ar[r]_f \ar[ru]^{\widetilde{f}}  & (M, b), 
}
$$
and which satisfies ${\mbox{\rm Im}}T_xf \subseteq \widetilde{f}(x)$ for any $x \in (N, a)$. Here ${\mbox{\rm Gr}}(n, TM)$ means the Grassmannian bundle over $M$ with fibers ${\mbox{\rm Gr}}(n, T_yM)$, 
the Grassmannians of $n$-dimensional subspaces in $T_yM, y\in M$. 
The condition on $\widetilde{f}$ is equivalent to 
that $\widetilde{f}$ is an integral mapping for the canonical distribution 
on the Grassmannian bundle. 

In this paper, in a dual manner to frontals, 
we introduce the notion of cofrontals: 
A map-germ $f : (N, a) \to (M, b)$ is called a {\it cofrontal} 
if $n \geq m$ and 
there exists an integrable vector-subbundle $K = K_f$ of $TN$ of corank $m$, 
which is regarded as a section 
$$
\xymatrix{
{\mbox{\rm Gr}}(n-m, TN) & { }
 \\
(N, a) \ar[r]_f \ar[u]^{K_f} 
& (M, b), 
}
$$
and which satisfies the condition 
$(K_f)_x \subseteq \Ker(T_xf)$ 
for any $x \in (N, a)$. 
We impose the integrability condition on $K$ in addition. 
If $f$ is fair (Definition \ref{Fair-cofrontal-germs}), i.e. the singular locus of the confrontal $f$ has no interior point nearby $a \in N$, then 
the integrability of the germ $K$ follows automatically. Moreover in this case $K$ is uniquely determined 
from the cofrontal $f$ (Lemma \ref{uniqueness-of-kernel-field}). 
In some sense, frontals are mappings such that the images of differentials 
are well-behaved, and cofrontals are mappings such that the kernels 
of differentials are well-behaved. 

A global mapping $f : N \to M$ is called a {\it frontal} (resp. a {\it cofrontal}) if all germs $f_a : (N, a) \to (M, f(a))$ 
of $f$ at every $a \in N$ are frontal (resp. cofrontal). Moreover a cofrontal $f$ is called {\it fair} 
if all germs $f_a$ at $f$ at every $a \in N$ are fair cofrontal germs (Definition \ref{Global-cofrontal-mappings}). 

Important examples of cofrontals are obtained as mappings which are constant along 
Seifert fibers (\cite{Audin}, cf. Example \ref{example-of-global-cofrontals}). 

We see that frontals and cofrontals are not stable except for the trivial cases, immersions and submersions and 
far from generic classes in the space of all $C^\infty$ mappings. 
Nevertheless we see that they enjoy rather interesting properties to be studied. 
For example, we see that 
any smooth map is approximated by a frontal or a cofrontal in $C^0$-topology, 
at least if the source manifold is compact 
(Proposition \ref{Topological-approximation}). 
In this paper we will describe such basic but interesting properties of cofrontals mainly. 

If $f : N \to M$ is a fair cofrontal, then 
the kernel field of $f$ exists uniquely and globally, and therefore 
the source manifold $N$ has a strict restriction if 
a global fair cofrontal exists on $N$. 
Note that, for given manifolds $N, M$ with $n < m$, if $N$ is compact then 
there exists a fair frontal $N \to M$ (Remark \ref{frontal-existence}). 

It is known that the local structures of fair (proper) frontals are understood by 
map-germs between spaces with the same dimension ($n = m$), 
together with the process of "openings" (\cite{Ishikawa18}\cite{Ishikawa18-2}). 
On the other hand the local structures of fair cofrontals turn to be reduced 
to the case $n = m$. In fact, 
as for the {\it source-local} problem, the classification of cofrontal singularities is reduced to the case 
$n = m$ completely 
(Proposition \ref{Criterion-of-cofrontality}, Lemma \ref{uniqueness-of-kernel-field}). 

Note that frontals were studied mainly in the case $m - n = 1$, i.e. the case of hypersurfaces, motivated by the study on wave-fronts (\cite{AGV}\cite{Ishikawa18}). 

In this paper, as for cofrontals, we study the cases of relative dimension $n - m = 1$.  
We provide a general {\it target-local} classification of fair cofrontals $N^{m+1} \to M^m$ 
with relative dimension $1$ under a mild condition. 
In fact the target-local classification problem of cofrontal mappings is reduced to that of 
the right-left classification of multi-germs $(\R^m, S) \to (\R^m, 0)$ 
together with a right symmetry of the multi-germ (Theorem \ref{classification-theorem}). 
It is interesting to apply the classification results 
of map-germs $(\R^m, 0) \to (\R^m, 0)$, in particular in the case $m = 2$ (see \cite{Whitney}\cite{Rieger}\cite{Saji}\cite{Kabata} for instance), to classifications of concrete classes of cofrontals. 

In \S \ref{cofrontals}, 
we introduce the notion of cofrontal map-germs comparing with that of frontals 
and clarify their local characters. In \S \ref{Global cofrontals}, 
we introduce global cofrontals and 
show some approximation result of mappings by frontals and cofrontals. 
After given several notions and examples related to fair frontals in 
\S \ref{Global fair cofrontals}, 
we give a classification of cofrontals of fiber-dimension one under the condition 
of \lq\lq reduction-finite" (Definition \ref{reduction-definition}, 
Definition \ref{reduction-finite-global}, 
Theorem \ref{classification-theorem}) 
in \S \ref{Classification of cofrontals of fiber-dimension one}. 

\

In this paper all manifolds and mappings are assumed to be smooth i.e. 
of class $C^\infty$ unless otherwise stated.

\section{Cofrontal singularities}
\label{cofrontals}

Let $N, M$ be smooth manifold of dimension $n$ and $m$ respectively, and 
$f : (N, a) \to (M, b)$ a smooth map-germ. Suppose $n \geq m$. 

\bef
\label{cofrontal-germ-definition}
{\rm (Cofrontal map-germ, kernel field.) \ 
The germ $f$ is called a {\it cofrontal map-germ} or a {\it cofrontal} in short,  if 
there exists a germ of smooth ($C^\infty$) integrable subbundle 
$K \subset TN$, $K = (K_x)_{x \in (N, a)}$,   
of rank $n-m$ such that 
$$
K_x \subseteq \ \Ker(T_xf : T_xN \to T_{f(x)}\R^m), 
$$
for any $x \in N$ nearby $a$. 
Here $T_xf : T_xN \to T_{f(x)}\R^m$ is the differential of $f$ at $x \in (N, a)$. 

Then $K$ is called a {\it kernel field} of the cofrontal $f$. 
}
\enf

Note that the kernel field is regarded as a section 
$K : (N, a) \to \Gr(n-m, TN)$ satisfying $(T_xf)(K_x) = \{ 0\}$, $x \in (N, a)$. 

Compare with the notion of frontals (cf. \cite{Ishikawa12}\cite{Ishikawa18}\cite{Ishikawa18-2}). 
Here we recall the definition of frontals: 
Let $f : (N, a) \to (\R^m, b)$ be a map-germ. Suppose $n \leq m$. 
Then $f$ is called a {\it frontal map-germ} or a {\it frontal} in short,  if 
there exists a smooth family of $n$-planes ${\widetilde f}(t) \subseteq T_{f(t)}\R^m$ 
along $f$, $t \in (N, a)$, 
satisfying the condition 
$
\Image(T_tf : T_tN \to T_{f(t)}\R^m) \subseteq \ {\widetilde f}(t) \ (\subset T_{f(t)}\R^m), 
$
for any $t \in (N, a)$. The family ${\widetilde f}(t)$ is called a {\it Legendre lift} of the frontal $f$. 

In some sense, cofrontals are dual objects to frontals. 

\bee
\label{examples-of-cofrontals}
{\rm 
(1) 
Any {\it immersion} is a frontal. 
The Legendre lift is given by $\widetilde{f} := (T_tf(T_tN))_{t \in (N, a)}$. 
Any {\it submersion} is a cofrontal.   
The kernel field $K$ is given by $K := (\Ker(T_xf))_{x \in (N, a)}$. 

(2) Any map-germ $(N, a) \to (M, b)$ between same dimensional manifolds $({\it n = m})$ is a frontal and a cofrontal simultaneously. 
In fact the Legendre lift is given by $\widetilde{f}(t) := T_{f(t)}M, t \in (N, a)$ and 
the kernel field $K$ is given by the zero-section of $TN$. 

(3) Any {\it constant map-germ} $(N, a) \to (M, b)$ 
is a frontal if $n \leq m$ and a cofrontal if $n \geq m$. In fact 
we can take any family of $n$-planes along the germ as a Legendre lift and 
any subbundle $K \subset TN$ of rank $n-m$ as a kernel field. 

}
\ene

\ber
{\rm 
As was mentioned, the differentials of cofrontals have a mild behavior. 
This reminds us Thom's $a_f$-condition: 
Let $f : N \to M$ be a smooth map, $X, Y$ submanifolds in $N$, and $x \in X \cap \overline{Y}$. 
Then $Y$ is $a_f$-regular over $X$ at $x$ if 
a sequence $y_i$ of points in $Y$ converges to $x$ and 
$$
\Ker(T_{y_i}(f\vert_Y)) \rightarrow T \subseteq T_xN, (i \to \infty), 
$$
then $\Ker(T_x(f\vert_X)) \subseteq T$. 

Let $f : N \to M$ is a cofrontal and take any fiber $X = f^{-1}(b), b \in M$. 
Then $X$ is a submanifold of $N$ and $Y = N \setminus Y$ is $a_f$-regular over $X$. 
}
\enr

%

%
%

\smallskip

Let ${\mathcal E}_{N,a} := \{ h : (N, a) \to \R\}$ 
denote the $\R$-algebra of smooth function-germs 
on $(N, a)$. 

Recall that {\it Jacobi ideal} $J_f$ of a map-germ $f : (N, a) \to (M, b)$ is defined as the ideal 
generated in ${\mathcal E}_{N,a}$ by all $\min \{n, m\}$-minor determinants of Jacobi matrix $J(f)$ of $f$. 
Note that $J_f$ is independent of the choices of local coordinates on $(N, a)$ and $(M, b)$. 

\bep
\label{Criterion-of-cofrontality}
{\rm (Criterion of cofrontality)} 
Let $f : (N, a) \to (M, b)$ be a map-germ with $n = \dim(N) \geq m = \dim(M)$. 
If $f$ is a cofrontal, then 
there exists a germ of submersion $\pi : (N, a) \to (\overline{N}, \overline{a})$ to an  
$m$-dimensional manifold $\overline{N}$ and a smooth map-germ 
$\overline{f} : (\overline{N}, \overline{a}) \to (M, b)$ such that $f = \overline{f}\circ\pi$. 
Moreover 
the Jacobi ideal $J_f$ of $f$ is principal, 
i.e. it is generated by one element. In fact $J_f$ is generated by $\lambda = \pi^*(\overline{\lambda})$ 
for the Jacobian determinant $\overline{\lambda}$ of $\overline{f}$. 

Conversely, if the Jacobi ideal $J_f$ is principal and the singular locus 
$$
S(f) = \{ x \in (N, a) \mid \rank(T_xf : T_xN \to T_{f(x)}M) < m \}
$$ 
of $f$ is nowhere dense in $(N, a)$, then $f$ is a cofrontal. 
\enp

\bef
\label{reduction-definition}
{\rm (Reductions of cofrontals.) \ 
We call $\overline{f}$ a {\it reduction} of the cofrontal-germ $f$. 
A germ of cofrontal $f : (N, a) \to (M, b)$ is called {\it reduction-finite} if 
a reduction $\overline{f} : (\overline{N}, \overline{a}) \to (M, b)$ of $f$ is ${\mathcal K}$-finite (or finite briefly), i.e. the dimension of 
$Q_{\overline{f}} := {\mathcal E}_{\overline{N}, \overline{a}}/\overline{f}^*(m_b)$ is finite, where 
$f^* : {\mathcal E}_{M,b} \to {\mathcal E}_{N, a}$ is the $\R$-algebra homomorphism defined 
by $f^*(h) = h\circ f$, 
and $m_b \subset {\mathcal E}_{M,b}$ is the maximal ideal of function-germs vanishing at $b$ (see \cite{Mather}\cite{GG}\cite{Wall}\cite{AGV}). 
}
\enf

\ber
\label{frontal}
{\rm 
In Lemma 2.3 of \cite{Ishikawa18-2}, 
it is shown that if $f : (N, a) \to (M, b), n \leq m$ is a frontal, 
then the Jacobi ideal $J_f$ is principal and that conversely 
if $J_f$ is principal and $S(f) = \{ x \in (N, a) \mid \rank(T_xf) < n \}$ is nowhere dense, 
then $f$ is a frontal. 
}
\enr

\ber
\label{finite-map-germs}
{\rm 
If $\overline{f} : (\R^m, 0) \to (\R^m, 0)$ is ${\mathcal K}$-finite, 
then the zero set of $\overline{f}$ is isolated and any nearby germ of $\overline{f}$ 
has the same property. 
The number of fibers of $\overline{f}$ is uniformly bounded by $\dim(Q_f)$ 
(Propositions 2.2, 2,4 of Ch.VII in \cite{GG}, see also \cite{EL}\cite{Khimshiashvili}). 
}
\enr

\noindent
{\it 
Proof of Proposition \ref{Criterion-of-cofrontality}}: \ 
Let $f$ be a cofrontal and $K$ be a kernel field of $f$. Since $K$ is integrable subbundle of $TN$ 
of rank $n - m$, there exists a 
submersion $\pi : (N, a) \to (\R^m, 0)$ such that $K_x = \Ker(\pi_* : T_xN \to T_{\pi(x)}\R^m$ 
for any $x \in (N, a)$, i.e. 
$\pi$-fibers form the foliation induced by $K$. 
Take any curve $\gamma : (\R, 0) \to N$ in a fiber of $\pi$. Then 
$(f\circ\gamma)'(t) = (T_{\gamma(t)}f)(\gamma \ '(t)) = 0$. Therefore $f$ is constant along the curve $\gamma$. 
Hence $f$ is constant on $\pi$-fibers. Then there exists a map-germ 
$\overline{f} : (\R^m, 0) \to (M, b)$ such that $f = \overline{f}\circ\pi$. 
Take a smooth section $s : (\R^m, 0) \to (N, a)$. Then 
$\overline{f} = \overline{f}\circ\pi\circ s = f\circ s$. Therefore $\overline{f}$ is a 
smooth map-germ. 

Take a system of local coordinates $x_1, \dots, x_m, x_{m+1}, \dots, x_n$ such that 
$\pi$ is given by 
$$
\pi(x_1, \dots, x_m, x_{m+1}, \dots, x_n) = (x_1, \dots, x_m)
$$ 
and therefore $K_x$ is generated by $\pa/\pa x_{m+1}, \dots, \pa/\pa x_{n}$ in $T_xN$. 
Then $f$ is expressed as 
$$
f(x_1, \dots, x_n) = (f_1(x_1, \dots, x_m), \dots, f_m(x_1, \dots, x_m)). 
$$
Then $J_f$ is generated by 
one element 
$\det(\pa f_i/\pa x_j)_{1 \leq i, j \leq m} 
= \pi^*(\det(\pa \overline{f}_i/\pa x_j)_{1 \leq i, j \leq m}) = \pi^*(\overline{\lambda})$ and therefore $J_f$ is a 
principal ideal in ${\mathcal E}_{N,a}$. 

Conversely suppose $J_f$ is a principal ideal generated by one element $\lambda \in J_f$ 
and $S(f)$ is nowhere dense. 
Denote by $\Gamma$ the set of subsets $I \subseteq \{ 1, 2, \dots, n\}$ 
with $\#(I) = m$. 
For a map-germ $f : (N, a) \to (M, b), n \geq m$ and $I \in \Gamma$, 
we set $D_I = \det(\pa f_i/\pa x_j)_{1 \leq i \leq m, j \in I}$ for some coordinates $x_1, \dots, x_n$ of $(N, a)$ 
and $y_1, \dots, y_m$ of $(M, b)$ with $f_i = y_i\circ f$. 
For any $I \in \Gamma$, there exists $h_I \in {\mathcal E}_a$ such that 
$D_I = k_I \lambda$. Since $S(f)$ is nowhere dense, 
there exists $I_0 \in \Gamma$ such that $D_{I_0} \not= 0$. 
Since $\lambda \in J_f$, 
there exists $\ell_I \in {\mathcal E}_a$ for any $I \in \Gamma$ 
such that $\lambda = \sum_{I \in \Gamma} \ell_I D_I$. Then 
$(1 - \sum_{I \in \Gamma} \ell_Ik_I)\lambda = 0$. 
If $k_I(a) = 0$ for any $I \in \Gamma$, then $1 - \sum_{I \in \Gamma} \ell_Ik_I$ is invertible in 
${\mathcal E}_a$, therefore $\lambda = 0$ and 
then we have $J_f = 0$. This contradicts to the assumption 
that $S(f)$ is nowhere dense. 
Hence there exists $I_0 \in \Gamma$ such that $(\ell_{I_0}k_{I_0})(a) \not= 0$. 
Then $k_{I_0}(a) \not= 0$. Therefore $J_f$ is generated by $D_{I_0}$. Hence 
$D_I = h_I D_{I_0}$ for any $I \in \Gamma$ with $h_{I_0}(a) = 1$. 
Then the Pl{\" u}cker-Grassmann 
coordinates $(h_I)_{I \in \Gamma}$ give a smooth section $K : (\R^a, a) \to 
\Gr(n - m, TN) \cong \Gr(m, T^*\R^a)$, which is regarded as a subbundle $K \subseteq TN$ 
of rank $n - m$ and $K_x \subseteq \Ker(T_xf)$ for any $x \in (N, a)$. 
Moreover $K_x$ coincides with $\Ker(T_xf)$ for $x \in (N \setminus S(f), a)$ 
and therefore $K$ is integrable outside of $S(f)$. Since $S(f)$ is nowhere dense, $K$ is integrable. 
Thus $f$ is a cofrontal map-germ with the kernel field $K$. 
\QED

\bec
Let $f : (N, a) \to (M, b)$ be a map-germ. Suppose $f$ is analytic and $J_f \not= 0$. 
Then $f$ is a frontal or a cofrontal if and only if $J_f$ is a principal ideal. 
\enc

\Proof
By Lemma \ref{Criterion-of-cofrontality} and Remark \ref{frontal}, 
if $f$ is a frontal or a cofrontal, then $J_f$ is principal. 
If $J_f$ is principal, $J_f \not= 0$ and $f$ is analytic, then $S(f)$ is nowhere dense. 
Thus by Lemma \ref{Criterion-of-cofrontality} and Remark \ref{frontal}, 
$f$ is a frontal if $n \leq m$ or a cofrontal if $n \geq m$. 
\QED

\bee
{\rm 
Let $f : (\R^3, 0) \to (\R^2, 0)$ be the map-germ given by 
$f(x_1, x_2, x_3) = (x_1^2 + x_2^2 + x_3^2, 0)$. Then $f$ is analytic and $J_f = 0$ is principal. 
However $f$ is not a cofrontal. 
In fact, suppose $f$ is a cofrontal and $K$ a kernel field of $f$ of rank $1$. 
Let 
$$
\xi(x) = \xi_1(x)\pa/\pa x_1 + \xi_2(x)\pa/\pa x_2 + \xi_3(x)\pa/\pa x_3, \ \xi(0) \not= 0,
$$ 
be a generator of $K$. 
Then $\xi_1(x)x_1 + \xi_2(x)x_2 + \xi_3(x)x_3$ is identically zero in a neighborhood of $0$ in $\R^3$. 
In particular we have $\xi_1(x_1, 0, 0)x_1 = 0$ and therefore $\xi_1(x_1, 0, 0) = 0$, so $\xi_1(0, 0, 0) = 0$. 
Similarly we have also $\xi_2(0, 0, 0) = 0$ and $\xi_3(0, 0, 0) = 0$. This leads a contradiction. 
}
\ene

\bef
\label{Jacobian}
{\rm (Jacobians of frontals and cofrontals.)\ 
Let $f : (N, a) \to (M, b)$ be a frontal or a cofrontal. 
Then a generator $\lambda \in {\mathcal E}_a$ 
of $J_f$ is called a {\it Jacobian} (or a {\it singularity identifier}) of the cofrontal $f$, 
which is uniquely determined from $f$ up to multiplication of a unit in ${\mathcal E}_a$. 
}
\enf

\ber
{\rm 
The singular locus $S(f) = \{ x \in (N, a) \mid \rank(T_xf : T_xN \to T_{f(x)}M) < \min\{ n, m\} \}$ 
of a frontal or a cofrontal $f$ is given by the zero-locus of the Jacobian $\lambda$ of $f$. 
}
\enr

\ber
{\rm
Let $f : (N, a) \to (M, b)$ be a cofrontal and $K$ a kernel field of $f$. Set 
$$
K^\perp_x := \{ \alpha \in T^*_xN \mid \alpha(v) = 0 \ {\mbox{\rm for any}} \ v \in K_x\}. 
$$
Then $K^\perp$ is a germ of subbundle of the cotangent bundle $T^*N$ of rank $m$. 
Let $\alpha_1, \alpha_2, \dots, \alpha_m$ be a local frame of $K^\perp$. 
Then there is a unique $\lambda \in {\mathcal E}_a$ such that
$$
df_1 \wedge df_2 \wedge \cdots \wedge df_m = \lambda \alpha_1\wedge \alpha_2 \wedge \cdots \wedge \alpha_m. 
$$
Then $\lambda$ generates $J_f$ and therefore $\lambda$ is a Jacobian of the cofrontal $f$. 
}
\enr

\bef 
\label{Fair-cofrontal-germs}
{\rm (Fair frontals and cofrontals.) \ }
{\rm 
A frontal or a cofrontal $f : (N, a) \to (M, b)$ is called {\it fair} if 
the singular locus $S(f)$ is nowhere dense in $(N, a)$. 
}
\enf

\ber
{\rm
A cofrontal $f$ is fair if and only if a reduction $\overline{f}$ (Definition \ref{reduction-definition}) is fair. 
In fact if $f =\overline{f}\circ\pi$ for a submersion-germ $\pi : (N, a) \to (\R^m, 0)$, we have 
$S(f) = \pi^{-1}(S(\overline{f}))$, and therefore $S(f)$ is nowhere dense in $(N, a)$ if and only if 
$S(\overline{f})$ is nowhere dense in $(\R^m, 0)$. 
If a cofrontal $f$ is reduction-finite (Definition \ref{reduction-definition}), 
then $f$ is fair, since a reduction $\overline{f}$ is 
${\mathcal K}$-finite so is necessarily fair. 
}
\enr

\ber
{\rm 
In \cite{Ishikawa18}\cite{Ishikawa18-2}, a frontal with nowhere dense singular locus was 
called {\it proper}. However in the global study the terminology \lq\lq proper" is rather confusing, 
in particular for the study of cofrontals, 
since its usage is different from the ordinary meaning of properness (inverse images of any compact is compact). 
Therefore in this paper we use the terminology \lq\lq fair" instead of \lq\lq proper". 
}
\enr

\bel
\label{uniqueness-of-kernel-field}
Let $f : (N, a) \to (M, b)$ be a fair cofrontal or $\dim(N) = \dim(M)$. 
Then the kernel filed $K$ of $f$ is uniquely determined and the reduction $\overline{f}$ of $f$ 
(Definition \ref{reduction-definition}) is uniquely determined up to right equivalence. 
\enl

\Proof
On the regular locus $N \setminus S(f)$, there is 
the unique kernel field $K$ defined by 
$K_x := \Ker(T_xf : T_xN \to T_{f(x)}M)$. 
Let $f$ be a fair cofrontal. Then $N \setminus S(f)$ 
is dense in $(N, a)$. Therefore the extension of $K$ to $(N, a)$ is unique if it exists. 
Let $n = m$. Then the unique unique kernel field $K$ is defined by the zero-section of $TN$ 
(Example \ref{examples-of-cofrontals} (2)). 
Then the submersion $\pi : (N, a) \to (\R^m, 0)$ induced by $K$ is uniquely determined 
up to left equivalence. Let $\pi' : (N, a) \to (\R^m, 0)$ be induced by $K$ and 
$\overline{f}$ and $\overline{f}'$ be both reductions of $f$ with 
$f = \overline{f} \circ \pi = \overline{f}'\circ \pi'$. 
Then $\pi' = \sigma\circ \pi$ for some 
diffeomorphism-germ $\sigma : (\R^m, 0) \to (\R^m, 0)$ and  
$\overline{f} \circ \pi = (\overline{f}'\circ\sigma)\circ\pi$. Since $\pi$ is a submersion, we have 
$\overline{f} = \overline{f}'\circ\sigma$. 
\QED

\

Let $f : (N, a) \to (M, b)$ be a cofrontal
(resp. a fair cofrontal) 
and $K : (N, a) \to {\mbox{\rm Gr}}(n-m, TM)$ be a kernel field of $f$. 
Recall that $K \subset TN$ is a germ of integrable subbundle of rank $n-m$. 

\bef
\label{adapted-coordinates}
{\rm (Adapted coordinates.)\ 
A system $(x_1, \dots, x_m, x_{m+1}, \dots, x_n)$ of local coordinates of $N$ centered at 
$a$ is called {\it adapted} to a kernel field $K$ of a cofrontal $f$, or simply, to $f$, 
if 
$$
\begin{array}{rcl}
K_x & = & \left\langle \left(\dfrac{\pa}{\pa x_{m+1}}\right)_x, 
\dots, \left(\dfrac{\pa}{\pa x_n}\right)_x\right\rangle_{\R}
= \{ v \in T_xN \mid dx_{1}(v) = 0, \dots, dx_{m}(v) = 0 \}, 
\end{array}
$$
for any $x \in (N, a)$. 
}
\enf

Since a kernel field $K$ of a cofrontal is assumed to be integrable, we have 

\bel
Any cofrontal $f : (N, a) \to (M, b)$ 
has an adapted system of local coordinates on $(N, a)$. 
\enl

\ber
{\rm 
For an adapted system of coordinates $(x_1, \dots, x_n, x_{n+1}, \dots, x_m)$ of $f$, 
the Jacobian $\lambda$ is given by the ordinary Jacobian $\frac{\pa(f_1, \dots, f_m)}{\pa(x_1, \dots, x_m)}$.
}
\enr

\section{Global cofrontals}
\label{Global cofrontals}

We will define the class of (co)frontal maps and fair (co)frontal maps, and we make clear the difference of these classes of mappings. 

\bef
\label{Global-cofrontal-mappings}
{\rm (Global cofrontal mappings.)\ 
Let $N, M$ be smooth manifolds of dimension $n, m$ respectively. 

Suppose $n \leq m$. 
A smooth mapping $f : N \to M$ is called a {\it cofrontal map} or a {\it cofrontal} briefly if 
the germ $f_a : (N, a) \to (M, f(a))$ at $a$ is a cofrontal for any $a \in N$ (Definition \ref{cofrontal-germ-definition}). 
A cofrontal $f : N \to M$ is called {\it fair} if $f_a$ is a fair cofrontal for any $a \in N$, i.e. if 
the singular locus $S(f) := \{ x \in N \mid \rank(T_xf : T_xN \to T_{f(x)}M) < m\}$ is nowhere dense in $N$ 
(Definition \ref{Fair-cofrontal-germs}). 

Suppose $n \geq m$. 
A smooth mapping $f : N \to M$ is called a {\it frontal map} or a {\it frontal} briefly if 
the germ $f_a : (N, a) \to (M, f(a))$ is a frontal for any $a \in N$. 
A frontal $f : N \to M$ is called {\it fair} if $f_a$ is a fair frontal for any $a \in N$, i.e. if 
the singular locus $S(f) := \{ x \in N \mid \rank(T_xf : T_xN \to T_{f(x)}M) < n\}$ is nowhere dense in $N$. 
}
\enf

\bee
\label{example-of-global-cofrontals}
{\rm
(1) Any submersion is a cofrontal. Any immersion is a frontal. 

(2) Any constant mapping $N \to M$ is a cofrontal of a frontal depending on $\dim(N) \geq \dim(M)$ or 
$\dim(N) \leq \dim(M)$. 

(3) Let ${\mathcal F})$ be 
a foliation of codimension $m$ on a manifold $N$ of dimension $n$. 
If a mapping $f : N^n \to M^m$ is constant on any leaf of ${\mathcal F}$, then $f$ is a cofrontal. 

(4) As a motivating example from symplectic geometry, consider a Lagrangian foliation ${\mathcal L}$ 
on a symplectic manifold $N^{2n}$ and a system of functions $f_1, \dots, f_n$ on $N$. Then 
$f = (f_1, \dots, f_n) : N \to \R^n$ is a cofrontal if $f$ is constant along each leaf of ${\mathcal L}$. 
}
\ene

First we observe \lq\lq unfair" (co)frontal maps are not so restrictive in topological or homotopical sense. 
In what follows we suppose $N$ is compact for simplicity. 

\bep
\label{Topological-approximation}
($C^0$-approximation.)\ 
Any proper smooth ($C^\infty$) map $f : N \to M$ is $C^0$-approximated by a frontal or a cofrontal $g : N \to M$. 
Any proper smooth map $f : N \to M$ is homotopic to a frontal or a cofrontal $g : N \to M$. 
\enp

\bee
{\rm
Let $S^2 \subset \R^3$ be the unit sphere and $g : S^2 \to \R$ 
the height function, i.e. $g(x_1, x_2, x_3) = x_3$. Then $g$ is never a cofrontal. 
Let $\varepsilon > 0$. 
Let $\varphi : [-1, 1] \to [-1, 1]$ be any smooth map satisfying that
$\varphi(y) = -1 (-1 \leq y -1 + \varepsilon), \varphi(y) = 1 (1 - \varepsilon \leq y \leq 1)$, 
and that $\varphi$ is a diffeomorphism 
from $(-1 + \varepsilon, 1 - \varepsilon)$ to $(-1, 1)$. 
Then $f = \varphi\circ g$ is a cofrontal. See the figure: In the right picture, 
$f$ restricted to the north (resp. south) gray part is constant. 

Note that $f$ can be taken to be arbitrarily near $g$ in $C^0$-topology. 

Similar construction can be applied to any proper Morse function $g : N \to \R$ and we have a cofrontal 
which is a $C^0$-approximation to $g$. 
}
\ene

\noindent
{\it Proof of Proposition \ref{Topological-approximation}}: 
Let $f : N \to M$ be a smooth mapping. 
Then 
$f$ is $C^\infty$-approximated by a mapping $f' : N \to M$ such that 
there exists a Thom stratification $({\mathcal S}, {\mathcal T})$ of $f'$. 
This is a consequence of the topological stability theorem (\cite{Mather70}\cite{Mather73}\cite{Mather76}\cite{GWPL}). 
In particular we have

(1) ${\mathcal S}$ is a Whitney stratification of $N$ 
and ${\mathcal T}$ is a Whitney stratification of $M$. 

(2) For any $S \in {\mathcal S}$, there exists a $T \in {\mathcal T}$ such that 
$f'\vert_S : S \to T$ is a surjective submersion. 

(3) The critical locus 
$$
\Sigma(f') := \{ x \in N \mid \rank(T_xf' : T_xN \to T_{f(x)}M) < m\}
$$
is a union of strata of ${\mathcal S}$ and, for any stratum $S \in {\mathcal S}$ in $\Sigma(S)$, 
$f\vert_S : S \to M$ is an immersion. 

(4) There exists 
a compatible tubular system $(\pi_S, \rho_S)_{S \in {\mathcal S}}$ for ${\mathcal S}$, i.e., 
a normed normal bundle $E_S \to S$ to $S$ in $N$, a positive smooth function 
$\varepsilon : S \to \R$ and a diffeomorphism $\phi_S : (E_S)_{<\varepsilon} \to N$ on on the image 
$U_S = \Phi_S((E_S)_{<\varepsilon})$, 
which is a tubular neighborhood of $S$ in $N$. Here 
$(E_S)_{<\varepsilon} := \{ v \in E_S \mid v \in T_xN, \Vert v\Vert < \varepsilon(x)\}$ 
and 
$\rho_S(v) = \Vert v\Vert^2$. 
The projection $\pi_S$ is regarded as the projection from the tubular neighborhood $U_S$ to $S$ 
via $\phi_S$, 
and $\rho_S$ is the squared norm function on a normed normal bundle of $S$, which is regarded 
as a function on the tubular neighborhood. 
Then compatibility condition means that, 
for any $S, S' \in {\mathcal S}$ with $S' \subseteq \overline{S}$, 
$$
\pi_{S'}\circ\pi_{S} = \pi_{S'}, \quad \rho_{S'}\circ\pi_{S} = \rho_{S}, 
$$
in the intersection of a tubular neighborhood of $S$ and that of $S'$, and that 
for any $S, S' \in {\mathcal S}$ with $\dim(S) = \dim(S')$, the intersection $U_S \cap U_{S'} = \emptyset$. 


For $0 \leq i \leq n$, denote by $S^{(i)}$ is the $i$-skeleton of ${\mathcal S}$, i.e. 
the union of all strata of ${\mathcal S}$ of dimension $\leq i$. 
Set $U^{(i)} = N \setminus S^{i}$, the union of all strata of ${\mathcal S}$ of dimension $> i$. 
We will modify $f'$ first on $U^{(n-1)}$ and then $U^{(n-2)}$ and so on to get the approximation $g$. 

Actually we perform as follows: First we suppose $n > m$. 
Let $\tau$ be a sufficiently small positive real number 
and $\tau = \tau_\delta : [0, 1] \to \R$ a smooth function such that $\tau(t) = 0 (0 \leq t < \delta), 
\tau(t) = 1 (1 - \delta < t \leq 1$. 
Define $f_{n-1} : U^{(n-1)} \to M$ by 
$$
f_{n-1}(x) := f\left(\phi_S\left(\tau\left(\frac{1}{\varepsilon(\pi_S(\phi_S^{-1}(x)))}\Vert\phi_S^{-1}(x)\Vert\right)\phi_S^{-1}(x)\right)\right), 
$$
for $x \in U_S \cap U^{(n-1)}$ with $\dim(S) = n-1$, and by $f_{n-1}(x) = f(x)$ otherwise. 
Note that, by the mapping $f_{n-1}$, the mapping $f$ is modified along 
points near $S \in {\mathcal S}$ with $\dim(S) = n-1$ contracts to $S$ 
and then mapped by $f$. The modified map $f_{n-1}$ is a smooth map and a cofrontal. 
Also we have that $f_{n-1}$ is homotopic to $f\vert_{U^{(n-1)}}$. 
Moreover note that $f_{n-1}$ is 
not a $C^\infty$-approximation 
but a $C^0$-approximation of $f$ on $U^{(n-1)}$. 
Define $f_{n-2} : U^{(n-2)} \to M$ by setting 
$f_{n-2}(x)$ similar as above 
for $x \in U_S \cap U^{(n-2)}$ with $\dim(S) = n-2$, by $f_{n-2}(x) = f_{n-1}(x)$ otherwise. 
Then $f_{n-2}$ is a smooth map, a cofrontal and a $C^0$ approximation of $f$ on $U^{(n-2)}$. 
Iterating this procedure we have $f_{0} : U^{(0)} \to M$ and finally $f_{-1} : U^{(-1)} = N \to M$, 
which is a smooth map, a cofrontal, a $C^0$ approximation of $f$ on $N$ and is homotopic to $f$. 

If $n = m$, then we have nothing to do. 

Suppose $n < m$. 
Note that in this case $\Sigma(f) = N$. 
Then by the same procedure as above, we have 
a {\it frontal} $f_{i}$ which is a $C^0$-approximation 
of $f\vert{U^{(i)}}$ and is homotopic to $f\vert{U^{(i)}}$ 
for $i = n-1, n-2, \dots, 0, -1$. 
Note that we may take as a Legendre lift of $f_{i}$ any extension of the Legendre lift 
of $f\vert_{S}, \dim(S) = i$ over $U_S$. 
Thus we have a frontal $f$ which is a 
$C^0$-approximation of $f$ and is homotopic to $f$. 
\QED

\section{Global fair cofrontals}
\label{Global fair cofrontals}

Contrary to the case of \lq\lq unfair" cofrontals, 
the following lemmata show that the sauce space of a {\it fair} cofrontal must be very restrictive. 

\bel
Let $f : N \to M$ be a fair cofrontal. Then there exists unique kernel field $K$ of $f$, 
i.e. there exists a unique integrable subbundle $K \subseteq TN$ of rank $n - m$ such that 
$K_x \subseteq \Ker(T_xf : T_xN \to T_{f(x)}M)$ for any $x \in N$. 
\enl

\Proof
By Lemma \ref{uniqueness-of-kernel-field}, the germ of kernel field is uniquely determined for each $x \in N$. 
By the local existence and {\it uniqueness}, we have the global existence of the kernel field of $f$. 
\QED

\bel
Let $f : N \to M$ be a fair cofrontal and $K$ the kernel field of $f$. 
Let ${\mathcal F}$ be the foliation induced by the integrable subbundle $K$ of $TN$ of rank $n - m$. 
Then the closure of any leaf of ${\mathcal F}$ is nowhere dense in $N$. 
\enl

\Proof
Let $L$ be a leaf of ${\mathcal F}$. Then $f$ restricted to $L$ is constant. (See the proof of 
Proposition \ref{Criterion-of-cofrontality}. Note that $L$ is assumed to be connected by the definition of leaves. )
Then $f$ restricted to the closure $\overline{L}$ of 
$L$ is constant just by the connectedness of $f$. Assume $\overline{L}$ has an interior point. 
Then also $S(f)$ necessarily has an interior point. This leads us to a contradiction with the fairness. 
\QED

\ber
\label{frontal-existence}
{\rm
Let $N, M$ be smooth manifolds with $\dim(N) = n \leq m = \dim(M)$. 
Suppose that $N$ is compact or both $N, M$ are non-compact. 
Then there exists a proper fair frontal $f : N \to M$. In fact take any closed submanifold $N' \subseteq M$ 
of dimension $n$ 
and its inclusion $i : N' \hookrightarrow M$. Take any proper smooth map $g : N \to N'$ whose singular locus 
$S(g)$ is nowhere dense, for instance, $g$ is a topologically stable map. Then 
$f := i\circ g$ is a proper fair frontal map. 
}
\enr

\

The following is clear. 

\bel
Let $g : L \to M$ be a cofrontal and $\pi : N \to L$ be a submersion. 
Then $g\circ \pi : N \to M$ is a cofrontal. $g\circ \pi$ is fair if and only if $g$ is fair. 
\enl

\bef
\label{irreducible-frontal}
{\rm (Reducible and irreducible cofrontals.)  \ 
Let $f : N \to M$ be a cofrontal with $\dim(N) = n > m = \dim(M)$. 
The frontal $f$ is called {\it reducible} if 
there exists a submersion $\pi : N \to \widetilde{N}$ to an $\ell$-dimensional manifold $\widetilde{N}$ 
and a cofrontal $g : \widetilde{N} \to M$ 
with $n > \ell \geq m$ such that $f = g\circ\pi$. 
A cofrontal is called {\it irreducible} if it is not reducible. 
}
\enf

\bep
Let $f : N \to M$ be a fair cofrontal with $n > m$ and $K$ its kernel field. 
If the leaf space form an $m$-dimensional 
manifold $\overline{N}$ and $\pi : N \to \overline{N}$ 
is a surjective smooth submersion such that $\Ker(\pi_*) = K$, 
then $f$ is reducible. In fact, there exists a smooth map 
$g : \overline{N} \to M$ such that $f = g\circ \pi$. 
\enp

\Proof
Since $f$ is constant on each leaf of $K$, we have a map $g : \overline{N} \to M$ 
such that $f = g\circ \pi$. Take any leaf $L$ of $K$ and any point $x \in L$, 
then the reduction of the germ of $f$ at $x$ is given by the germ of $g$ at $L 
\in \overline{N}$. Therefore $g$ is smooth at $L$ (see 
Proposition \ref{Criterion-of-cofrontality}). 
Thus $g$ is a smooth map. 
\QED

\bee
{\rm
(Irreducible fair frontals.) \ 
(1)
Let the open M{\" o}bius band $N$ is given as the quotient of $\R^2$ by the cyclic action generated by the 
transformation 
$\psi : \R^2 \to \R^2$, $\psi(x_1, x_2) = (x_1 + 1, -x_2)$. 
Then $(x_1, x_2) \mapsto x_2^2$ induces a well-defined map $f : N \to \R$ which is 
an irreducible fair cofrontal. 

(2)
Let $T = \R^2/\ZZ^2$ be the torus. Let $K$ be the Klein bottle defined as the quotient by the involution 
$\varphi : T \to T$, $\varphi([x_1, x_2]) := [x_1 + \frac{1}{2}, 1 - x_2]$. 
Define $f : K \to \R$ by $f([[x_1, x_2]]) := (x_2 - \frac{1}{2})^2$. Here $[x_1, x_2]$ (resp. 
$[[x_1, x_2]]$) be the point on $T$ (resp. $K$) represented by $(x_1, x_2) \in \R^2$. 
Then $f$ is well-defined smooth mapping which is an irreducible fair cofrontal. 
}
\ene

\bee
{\rm
(Cofrontal of reduction-non-finite.)
Let $\varphi : \R \to \R$ be a smooth such that $\varphi(t) = t, (t < \frac{1}{3}, t\frac{2}{3} < t)$,  
$\frac{1}{3} < \vert \varphi(t)\vert < \frac{2}{3}$ and that $\varphi^{-1}(\frac{1}{2})$ is an 
infinite set having just one point $t = \frac{1}{2}$ as a non-isolated point. Then 
define $f : T^2 = \R^2/\ZZ^2 \to S^1 = \R/\ZZ$ by $f([t_1, t_2]) = [\varphi(t_2)]$ modulo $\ZZ$, where $t_2 \in [0, 1]$. 
Then $f$ is a fair cofrontal such that $f$ is not reduction-finite and 
the fiber $f^{-1}([\frac{1}{2}])$ has infinite many connected components. 
}
\ene


\section{Classification of cofrontals of fiber-dimension one} 
\label{Classification of cofrontals of fiber-dimension one}

To give a target-local classification theorem of cofrontals, we start with an algebraic consideration. 
Let us denote by ${\mbox{\rm Diff}}(N, a)$ the group of diffeomorphisms $(N, a) \to (N, a)$. 

\bef
{\rm 
Let $f : (N, a) \to (M, b)$ be a smooth map-germ. Then the {\it right symmetry group} $G_f$ of $f$ 
is defined by 
$$
G_f := \{ \sigma \in {\mbox{\rm Diff}}(N, a) \mid f\circ\sigma = f \}. 
$$
}
\enf

\bel
Let $f : (N, a) \to (M, b)$ and $g : (N', a') \to (M', b')$ be smooth map-germs. 
If $f$ and $g$ are right-left equivalent (${\mathcal A}$-equivalent), then 
$G_f$ and $G_g$ are isomorphic as groups. 
\enl

\Proof
Suppose $\tau\circ = g\circ\sigma$ for diffeomorphism-germs $\sigma : (N, a) \to (N', a')$ and 
$\tau : (M, b) \to (M', b')$. Let $\varphi \in G_f$. 
Then 
$$
g\circ(\sigma\circ\varphi\circ\sigma^{-1}) = (g\circ\sigma)\circ\varphi\circ\sigma^{-1} 
= (\tau\circ f)\circ\varphi\circ\sigma^{-1} = \tau\circ (f\circ\varphi)\circ\sigma^{-1}
\tau\circ f\circ \sigma^{-1} = g. 
$$
Therefore $\sigma\circ\varphi\circ\sigma^{-1} \in G_g$. The correspondence 
$G_f \to G_g$ defined by $\varphi \mapsto \sigma\circ\varphi\circ\sigma^{-1}$ induces 
a group isomorphism. 
\QED

\bee
{\rm
(1) Let $f : (\R^2, 0) \to (\R^2, 0)$ be a {\it fold} which is defined by $f(x_1, x_2) = (x_1, x_2^2)$. 
Then the right symmetry group $G_f \cong \ZZ/2\ZZ$. 

(2) Let $f : (\R^2, 0) \to (\R^2, 0)$ be a {\it cusp} which is defined by $f(x_1, x_2) = (x_1, x_2^3 + x_1x_2)$ (\cite{Whitney}). 
Then $G_f$ is trivial, i.e. $G_f$ consists of only the identity map-germ on $(\R^2, 0)$. 

(3) Let $f : (\R^2, 0) \to (\R^2, 0)$ be defined by $f(x_1, x_2) = (x_1^2, x_2^2)$. Then 
$G_f \cong \ZZ/2\ZZ \times \ZZ/2\ZZ$. 

(4) Let $f : (\R^2, 0) = (\C, 0) \to (\C, 0) = (\R^2, 0)$ be defined by $f(z) = z^\ell, (z \in \C)$. 
Then we have $G_f \cong \ZZ/\ell\ZZ$. 

(5) 
Let $G$ be a finite reflection group on $\R^n$ and $h_1, h_2, \dots, h_n$ be a system of generators  
of the invariant ring of $G$ consisting of homogeneous polynomials (cf. Chevalley's theorem \cite{Humphreys}). 
Then the right symmetry group $G_h$ of $h = (h_1, \dots, h_n) : (\R^n, 0) \to (\R^n, 0)$ 
is isomorphic to $G$. 
}
\ene

\bep
\label{construction-theorem}
{\rm (Construction of cofrontals of fiber-dimension one) \ } 
Let $h : (\R^m, 0) \to (M, b)$ be a smooth map-germ 
and $\sigma \in G_h$. 
Let $h : U \to M$ and $\sigma : U \to U$ be representatives of $h$ and $\sigma$ respectively such that $h\circ\sigma = h$ on $U$. 
Set $N = ([0, 1]\times U)/\!\sim$ where $(0, \overline{x}) \sim (1, \sigma(\overline{x}))$. 
Then $N$ is a $(m+1)$-dimensional manifold and 
$f = f_{h, \sigma} : N \to M, f([t, \overline{x}]) = h(\overline{x})$ is well-defined and is a  cofrontal. 

In general, let $h_1, \dots, h_s : (\R^m, 0) \to (M, b)$ 
be smooth map-germs and $\sigma_i \in G_{h_i}, (1 \leq i \leq s)$. 
Let $h_i : U_i \to M$ and $\sigma_i : U_i \to U_i$ be representatives of $h_i$ and $\sigma_i$, $(1 \leq i \leq s)$ 
respectively such that $h_i\circ\sigma_i = h_i$ on $U_i$. 
Set $N_i = ([0, 1]\times U_i)/\!\sim$ where $(0, \overline{x}) \sim (1, \sigma_i(\overline{x}))$. 
Take the disjoint union $N = \bigcup_{i=1}^s N_i$. 
which is an $(m+1)$-dimensional manifold. 
Define $f = f_{h_1, \dots, h_s; \sigma_1, \dots, \sigma_s}  : N \to M$ by 
$f([t, \overline{x}]) = h_i(\overline{x})$ for $[t, \overline{x}] 
\in N_i$, $1 \leq i \leq s$. 
Then $f$ is well-defined and $f$ is a cofrontal. 
\enp


\Proof
Since $\sigma_i \in G_{h_i}$, $f$ is well-defined and smooth. Moreover the $t$-direction defines 
well-defined subbundle $K \subset TN$ of ranks $1$. Since $f$ is constant along $K$, we see 
that $f$ is a cofrontal. 
\QED

\ber
{\rm
The cofrontal $f$ in Proposition \ref{construction-theorem} is fair if and only if 
all $h_i, i = 1, \dots, s$ are fair. Moreover $f$ is reduction-finite if and only if all $h_i, i = 1, \dots, s$ 
are ${\mathcal K}$-finite. 
}
\enr

\bef
\label{reduction-finite-global}
{\rm
(Reduction-finite cofrontals.)\ 
A cofrontal $f : N \to M$ is called {\it reduction-finite} if any germ 
of cofrontal $f_a : (N, a) \to (M, f(a))$ is reduction-finite 
in the sense of Definition \ref{reduction-definition}. 
}
\enf

\bef
{\rm (\cite{Saeki}) \ 
Let $f : N \to M, f' : N' \to M'$ be smooth map-germs and $b \in M, b' \in M'$. 
Then the germ of $f$ over $b$ is right-left equivalent the germ of $f'$ over $b'$ if 
there exists an open neighborhood $U$ of $b$ in $N$, an open neighborhood $U'$ of $b'$, 
a diffeomorphism $\Phi : f^{-1}(U) \to f'^{-1}(U')$ and a diffeomorphism $\varphi : U \to U'$ 
such that the diagram 
$$
\xymatrix{
f^{-1}(U) \ar[d]_{f} \ar[r]^{\Phi} &  f'^{-1}(U') \ar[d]^{f'}
 \\
U \ar[r]_{\varphi}  & U'
}
$$
commutes. 
}
\enf

\ber
{\rm 
The right-left equivalence class of the germ of $f_{g, \sigma} : N \to M$ over $b \in M$ 
in Proposition \ref{construction-theorem} depends only on 
the right-left equivalence class of the germ $g$ and 
the conjugacy class of $\sigma$ in $G_g$. 
Similarly the 
right-left equivalence class of the germ of 
$f_{h_1, \dots, h_s; \sigma_1, \dots, \sigma_s} : N \to M$ over $b \in M$ 
the right-left equivalence class of the multi-germ $(g_1, \dots, g_s)$ from 
the disjoint union of $s$-copies of $(\R^m, 0)$ to $(M, b)$,  
and the conjugacy classes of $\sigma_i$ in $G_{g_i}$. 
}
\enr

\bet
\label{classification-theorem}
{\rm (Classification theorem of cofrontals with one-dimensional fibers.)}
Let $N$ be a compact smooth manifold of dimension $m+1$, 
and $M$ a smooth manifold of dimension $m$. 
Let $f : N \to M$ be any reduction-finite cofrontal and $b \in M$. 
Then the germ $f$ over $b$ is right-left equivalent to the germ 
$f_{h_1, \dots, h_s; \sigma_1, \dots, \sigma_s}$ over $b$ for some 
non-negative integer $s$, ${\mathcal K}$-finite map-germs 
$h_i : (\R^m, 0) \to (M, b)$ and elements $\sigma_i \in G_{g_i}$ 
of finite order 
$(1 \leq i \leq s)$. 
\ent

\bel
\label{representatives-of-finite-germ}
Let $h : (\R^m, S) \to (M, b)$ be a multi-germ with $S = \{ x_1, \dots, x_s\}$. 
Suppose all germ $h_i = h_{x_i} : (\R^m, x_i) \to (\R^m, 0)$ are ${\mathcal K}$-finite. 
Let $\sigma_i \in G_{h_i}$. Then there exist open neighborhood $V$ of $b$, 
open neighborhood $U_i$ of $x_i$ and representatives $h_i : U_i \to M$ of $h_i$ and 
$\sigma_i : U_i \to U_i$ such that $h_i^{-1}(V) = U_i$ and $h_i\sigma_i = h_i$ on $U_i$ for $i = 1, \dots, s$. 
\enl

\Proof
Let $g : (\R^m, 0) \to (M, b)$ be a ${\mathcal K}$-finite map-germ. Then $m_0^k \subset f^*(m_b)$ 
for some positive integer $k$. Then there exist $\alpha > 0,  C > 0$ such that $C\Vert x\Vert^\alpha \leq 
\Vert g(x)\Vert$ on a neighborhood of $0 \in \R^m$. Therefore, for any representative $g : W \to M$ of 
$g$, and for any neighborhood $W'$ of $0$ with $W' \subseteq W$, there exists an open neighborhood 
$V$ of $b$ such that $g^{-1}(V) \subset W'$. For each $h_i$ take such an open neighborhood $V_i$ of $b$ 
such that $h_i\circ \sigma_i = h_i$ holds on $h^{-1}(V_i)$. 
Then set $V = \cap_{i=1}^s V_i$ and $U_i = h_i^{-1}(V)$. Then $\sigma(U_i) = U_i$ and 
$h_i\circ\sigma_i = h_i$ on $U_i$. 
\QED

\bel
\label{finite-fiber}
Let $f : N \to M$ be a cofrontal of reduction-finite, $N$ compact and $b \in M$. 
Then the fiber $f^{-1}(b)$ over $b$ consists of a finite number of disjoint circles in $N$. 
Each connected component has an open neighborhood consist of leaves of a kernel field of $f$. 
\enl

\Proof
First we remark that, since the cofrontal $f$ is reduction-finite, 
$f$ is fair and therefore there exists the unique global kernel field $K$ of $f$. 
Take any $a \in N$. Let $L$ be the leaf through $a$ of the foliation ${\mathcal F}$ defined by $K$. 
Then the germ $f_a$ has a ${\mathcal K}$-finite reduction. Then there exists an adapted 
open neighborhood $W_a$ of form $U \times V, \dim(U) = m, \dim(V) = n - m$ 
such that $\{p\}\times V$ is contained in a leaf of ${\mathcal F}$ for any $p \in U$. 
If we take $U$ sufficiently small, then $W_a \cap f^{-1}(b) = W_a \cap L$, since $f$ is reduction-finite. 
Set $W = \cup_{a \in L} W_a$. Then $W$ is an open set in $N$ and $W \cap f^{-1}(b) = W \cap L$. 
Thus we have seen that $L$ is a closed, therefore compact submanifold in $N$. In particular 
$L$ is diffeomorphic to the circle $S^1$. Moreover $L$ has an open neighborhood 
consists of leaves of ${\mathcal F}$. Then we have that the number of connected components is finite. 
\QED

\

\noindent
{\it Proof of Theorem \ref{classification-theorem}}: 
Since $N$ is compact and the cofrontal 
$f$ is reduction-finite, $f^{-1}(b)$ consists of a finite number of disjoint circles $L_1, L_2, \dots L_s$ 
in $N$ by Lemma \ref{finite-fiber}. 
Each $L_i$ has an open neighborhood $W_i$ consisting of leaves of the foliation ${\mathcal F}$ 
of the kernel field $K$ of $f$. By taking each $W_i$ small enough, we have that 
$W_i \cap W_j = \emptyset$ for $i \not= j$. Now $f$ is constant on each leaf of ${\mathcal F}$. 
Take a transversal $U_i$ of dimension $m$ to the leaves on $W_i$ through a point on $L_i$. 
Then we have the Poincar{\' e} map $\sigma_i : U_i \to U_i$ by moving along leaves of ${\mathcal F}$. 
We have that $f\circ\sigma_i = \sigma_i$ on $U_i$. Set $h_i = f\vert_{U_i}$. 
Then we have that 
the germ of $f$ over $b$ is right-left equivalent to $f_{g_1, \dots, g_s; \sigma_1, \dots, \sigma_s}$ over $b$. 
Taking $U_i$ sufficiently small, then the number of fibers of $h_i$ is bounded (Remark \ref{finite-map-germs}). 
Then any $\sigma_i$-orbit on $U_i$ has bounded period. Therefore 
$\sigma_i$ must be of finite order. 
\QED

\end{document}